\documentclass{article}
\usepackage{latexsym,amsmath,amsthm,amssymb}

\newcommand{\dR}{\ensuremath{\mathbb{R}}} 
\newcommand{\R}{\dR}

\newcommand{\var}{\mathbf{Var}} 
\newcommand{\osc}{\mathbf{Osc}} 
\newcommand{\capa}{\mathrm{Cap}} 

\newtheorem{theorem}{Theorem}
\newtheorem{proposition}[theorem]{Proposition}

\newtheorem{corollary}[theorem]{Corollary}

\theoremstyle{definition} 
\newtheorem*{assumption}{Hypothesis}

\theoremstyle{remark} \newtheorem{remark}{Remark}
\newtheorem{example}[remark]{Example}

\begin{document}   
\title{Concentration for independent random variables \\  with heavy tails}
\author{F. Barthe, P. Cattiaux and C. Roberto}
\maketitle

\begin{abstract}
   If a random variable is not exponentially integrable, it is known 
that no concentration inequality holds for an infinite sequence of 
independent copies. Under mild conditions, we establish concentration
inequalities for finite sequences of $n$ independent copies, with good
dependence in $n$. 
\end{abstract}

\section{Introduction}
This paper continues the study of the concentration of measure phenomenon
for product probability measures. A detailed account of this topic and its 
applications is given in \cite{ledoCMP}. Let us recall an important method
for this problem: if $\mu$ (say on $\dR^{d}$) satisfies a spectral gap  (or
Poincar\'e) inequality
$$ \var_\mu(f)\le C\int|\nabla f|^{2} d\mu,
\quad \mathrm{for\;  all\; locally\;  Lipschitz\; } f:\dR^{d}\to\dR,$$
then Lipschitz functions are exponentially concentrated \cite{gromm83taii,borou83icnd}.
 More precisely
every $1$-Lipschitz function $F$ in the Euclidean distance, with median $m_F$,
 satisfies
$\mu(|F-m_F|>t)\le 6 \exp(-t/(2\sqrt{C}))$ for $t>0$. Since the Poincar\'e inequality
has the so-called tensorisation property, the same property holds for 
$\mu^{n}$ for all $n\ge 1$. Similarly, the logarithmic Sobolev inequality
(see e.g. \cite{ledoCMLS}) yields dimension free Gaussian concentration,
 whereas recent inequalities devised by Lata{\l}a and Oleszkiewicz \cite{latao00bsp}
provide intermediate rates, see also \cite{bartr03sipm,wang03gbti,bartcr04iibe}.
 Note that these results only concern distributions
with exponential or faster decay. This was explained by Talagrand \cite{tala91niic}.
 Together with his famous result for products of exponential laws
 he observed the following: if $\mu$ is a probability measure on $\dR$ such that
there exist $h>0$ and $\varepsilon_{1/2}>0$ such that for all $n\ge 1$ and
all $A\subset \dR^{n}$ with $\mu^{n}(A)\ge \frac12$, one has
$$ \mu^{n}(A+[-h,h]^{n})\ge\frac12+\varepsilon_{1/2}$$
then $\mu$ has exponential tails, that is there exist positive constants $C_1,C_2$
such that $\mu([x,+\infty))\le C_1 e^{-C_2 x}$,
$x\in\dR$. A similar property for all $p\in(0,1)$ instead of just $p=1/2$ implies
that $\mu$ is the image of the symmetric exponential law by a map with
finite modulus of continuity, as Bobkov and Houdr\'e proved \cite{bobkh00wdfc}.

Thus when the tails of $\mu$ do not decay exponentially fast, there is no hope 
for dimension free concentration. This paper provides positive results in this 
case by investigating the size of enlargement $h_n$ necessary to ensure
a rise of the measure in dimension $n$. 
We study the more natural and also more difficult 
notion of Euclidean enlargement, and estimate $h_n$ such that 
$\mu^{n}(A)\ge 1/2$ implies $\mu^{n}(A+h_n B_2^{n})\ge \frac12+\varepsilon$,
where $B_2^{n}$ is the Euclidean unit ball. By the above results we know that
$h_n$ has to tend to infinity as the dimension $n$ increases. 
This question can be reformulated in terms of functions: we are looking for $h_n$
such that for all $n$ and all $1$-Lipschitz functions $F:\dR^{n}\to \dR$ with
median $m_F$, one has $\mu^{n}(F-m_F >h_n)\le \frac12-\varepsilon$.

\medskip
We work in the setting of a Riemannian manifold $(M,g)$ with a Borel probability 
measure which is absolutely continuous with respect to the volume measure. 
Our approach is based on the weak spectral gap inequality introduced by 
R\"ockner and Wang \cite{rockw01wpil}. In this remarkable paper, these authors 
 provide several 
necessary conditions for a measure to satisfy such a property,  consequences for
the corresponding semi-group and isoperimetric inequalities (see also
\cite{aida01egss,wangz03wpid} for other developments). Our results complete
and sharpen some of theirs. In Section 2 we give a characterization 
of measures on the real line with a weak spectral gap inequality.
Section~3 shows that this functional inequality has a defective tensorisation
property. We deduce isoperimetric and concentration inequalities for products in
 Sections~4 and 5. We illustrate our results with the examples of the power laws 
$\alpha (1+|t|)^{-1-\alpha}dt/2$ for $\alpha>0$ and the exponential type laws
$\exp(-|t|^{p})dt/(2\Gamma(1+1/p))$ for $p\in (0,1)$. The latter should be of importance
in the study of $p$-convex sets, as their analogues for $p\ge 1$ were in convex 
geometry (see e.g. \cite{schez00clb}). We discuss our concentration
consequences of the weak Poincar\'e inequality, in comparison with the ones
of the recent article \cite{wangz03wpid}. Our results are stronger, but the 
argument of Wang and Zhang can be improved in order to recover ours, and
actually a slightly better though less explicit bound. 
The final section illustrates our method on a wide family of measures extending
the laws $c_p \exp(-|t|^{p})dt$, $p\in(0,1)$.

\medskip
Let $\mu$ be an absolutely continuous probability measure on a
Riemannian manifold $M$. The modulus of gradient of a locally Lipschitz function 
$f:M\to \mathbb R$ can be defined as a whole by 
$$ |\nabla f|(x)=\limsup_{y\to x}\frac{|f(x)-f(y)|}{d(x,y)}$$
where $d$ is the geodesic distance.
Following R\"ockner and Wang, we say  $\mu$ satisfies a weak Poincar\'e 
inequality if there exists a function $\beta:(0,+\infty)\to\dR^{+}$ such that 
 every locally Lipschitz function  $f:M\to \dR$ satisfies for all $s>0$
the inequality
$$ \var_\mu(f) \le \beta(s) \int |\nabla f|^{2}d\mu + s \,\osc(f)^{2}.$$
Here $\osc(f)=\sup f -\inf f$ is the total oscillation of the function $f$.
The above mentioned authors used instead the quantity $\| f-\int f\,d\mu\|_\infty$.
When this $L_\infty$ essential supremum norm is with respect to the volume measure, the
two quantities are the same up to a factor 2. We shall assume as we may, that $\beta$
is non-increasing. Since $\var_\mu(f) \le \osc(f)^{2}/4$, the inequality is 
  trivial  when $s\ge 1/4$. In other words, one may set $\beta(s)=0$ for $s\ge 1/4$.
 The real content of the inequality is when $s$ is close to $0$.
If $\lim_{s\to 0}\beta(s) =b$, $b>0$ 
then the measure satisfies a classical Poincar\'e or 
spectral gap inequality. Otherwise the speed of convergence to $+\infty$ is of 
great  interest.

\section{A measure-capacity criterion}
This section provides an equivalent form of the weak Poincar\'e inequality,
in terms of a comparison between capacity of sets and their measure.
This point of view was put forward in \cite{bartcr04iibe} in order to
give a natural unified presentation of the many functional inequalities
appearing in the field.  In 
dimension 1 this leads to a very effective necessary and sufficient condition
for a measure to satisfy such an inequality, with a precise estimate of the 
function $\beta$. This completes the work by R\"ockner and Wang where several
necessary conditions were provided.

In the following, $\mathbf 1_S$ denotes the characteristic function of a set $S$,
and $f_{|S}$ is the restriction of the function $f$ to the set $S$.
Given measurable sets $A\subset \Omega$, the capacity $\capa_\mu(A,\Omega)$,
is defined as 
\begin{eqnarray*}
\capa_\mu(A,\Omega) &=& 
    \inf\left\{\int|\nabla f|^{2}d\mu;\; f_{|A}\ge 1,\; f_{|\Omega^{c}}=0
           \right\} \\
  &=&\inf\left\{\int|\nabla f|^{2}d\mu;\; 
         \mathbf 1_A \le f \le \mathbf 1_{\Omega} \right\},
\end{eqnarray*}
where the infimum is over locally Lipschitz functions. The latter equality
follows from an easy truncation argument, reducing to functions with values
in   $[0,1]$. Finally we defined in \cite{bartr03sipm} the capacity of $A$ with respect to 
$\mu$ when $\mu(A) < 1/2$ as 
$$ \capa_\mu(A):=\inf\{ \capa(A,\Omega); \; A\subset \Omega,\; \mu(\Omega)\le 1/2\}.$$

\begin{theorem}
Assume that for every $f:M\to \dR$ and every $s\in (0,1/4)$ one has
$$ \var_\mu(f) \le \beta(s) \int |\nabla f|^{2}d\mu + s \,\osc(f)^{2}.$$
Then for every measurable $A\subset M$ with $\mu(A)<1/2$, one has
$$ \capa_\mu(A)\ge \frac{\mu(A)}{4\beta(\mu(A)/4)}\cdot $$
\end{theorem}

\begin{proof}
   We start with assuming the weak Poincar\'e inequality. Let $A\subset \Omega$,
where $\mu(\Omega)\le 1/2$. Let $f$ be a locally Lipschitz function satisfying
$ \mathbf 1_A \le f \le \mathbf 1_{\Omega}$. By Cauchy-Schwarz inequality,
$$ \left( \int f\, d\mu \right)^{2}=\left( \int f \mathbf 1_\Omega\, d\mu \right)^{2}
 \le \mu(\Omega) \int f^{2}d\mu.$$
Therefore $\var_\mu(f)\ge \mu(\Omega^{c})\int f^{2 } d\mu\ge \int f^{2}d\mu/2.$
Since the oscillation of $f$ is at most $1$, the weak Poincar\'e inequality
yields for $s\in (0,1/4)$
$$ \frac12 \mu(A) \le \frac12 \int f^{2} d\mu 
\le \beta(s) \int |\nabla f|^{2}d\mu + s.$$
This is valid for arbitrary $f$ with $ \mathbf 1_A \le f \le \mathbf 1_{\Omega}$.
Hence we get  $$ \frac12 \mu(A) \le \beta(s)\, \capa_\mu(A,\Omega)+s.$$
Taking the infimum over sets $\Omega$ with measure at most $1/2$ and containing $A$,
we obtain for any $s\in(0,1/4)$
$$ \frac{1}{\beta(s)}\left(\frac{\mu(A)}{2}-s \right)_+ \le \capa_\mu(A).$$
Note that as a function of $\mu(A)$ the above lower bound vanishes before $2s$
and then increases with slope $1/(2\beta(s))$. Taking supremum over $s$ 
yields general lower bounds of the capacity by convex functions of the measure,
vanishing at $0$. More precisely we arrived at $\capa_\mu(A)\ge \tilde{\beta}(\mu(A))$,
where for $a\in(0,1/2)$, $$\tilde{\beta}(a)
=\sup_{s\in(0,1/4)} \left(\frac{ \frac{a}{2}-s}{\beta(s)} \right)_+ = 
\sup_{s\in(0,a/2)} \frac{ \frac{a}{2}-s}{\beta(s)}  \cdot$$
Note that 
$$ \frac{a}{4\beta(a/4)} \le \tilde{\beta}(a) \le \frac{a}{2\beta(a/2)},  $$
where the lower bound corresponds to the choice $s=a/4$ and the upper bound
relies on the non-increasing property of $\beta$. When this function satisfies
a doubling condition ($\beta(2x)\ge c \beta(x)$) then the above bounds 
are the same up to a multiplicative constant.
\end{proof}

\begin{theorem}\label{crit-capa}
   Assume that $\gamma$ is a non-increasing positive  function on $(0,1/2)$.
 If every measurable $A\subset M$ with $\mu(A)\le 1/2$ verifies
 $$ \capa_\mu(A)\ge \frac{\mu(A)}{\gamma(\mu(A))},$$
then for every locally Lipschitz function $f$  and every $s\in(0,1/4)$ one has
$$ \var_\mu(f) \le 12\gamma(s) \int |\nabla f|^{2}d\mu+ s\, \osc(f)^{2}.$$
\end{theorem}
\begin{proof}
 Fix $s\le 1/4$. Let $m$ be a median of $f$ under $\mu$. Denote
$\Omega_+=\{f>m\}$ and $\Omega_-=\{f<m\}$. Then
$$\var_\mu(f) \le \int (f-m)^{2}d\mu =  \int_{\Omega_+} (f-m)^{2}d\mu+
                                           \int_{\Omega_-} (f-m)^{2}d\mu.$$
We work separately on each of the latter two integrals.
Consider $g=(f-m)_+$ as a function defined on $\Omega_+$.
 Let $c=\inf\{t\ge 0;\; \mu(g^{2}>t)\le s \}$.
If $c=0$ then $\mu(g>0)\le s$ and $\int_{\Omega_+} g^{2}d\mu \le s \max g^{2}$
and we are done for this half of space. 
Otherwise $\mu(g^{2}>c)\le s$ and $\mu(g^{2}\ge c)\ge s$. 
By our structural hypothesis of a Riemannian manifold with an absolutely 
continuous measure we can find a set $\Omega_0$ with
$\{g^{2}>c\} \subset \Omega_0 \subset \{g^{2}\ge c \}$ and $\mu(\Omega_0)=s$.
Let $\rho>1$. For $k<0$ and integer, 
define $\Omega_k=\{g^{2}\ge c\rho^{k}\}$. Then
\begin{eqnarray*}
   \int_{\Omega_+}g^{2}d\mu &=& \int_{\Omega_0}g^{2}d\mu+\sum_{k<0}
           \int_{\Omega_k\setminus \Omega_{k+1}}g^{2}d\mu \\
  &\le & s \,\sup (f-m)_+^{2}+
  \sum_{k<0} c \rho^{k+1} \Big(\mu(\Omega_k)-\mu(\Omega_{k+1})\Big)
\end{eqnarray*}
The second term is dealt with by Abel summation:
\[ \sum_{k<0}  \rho^{k+1} \Big(\mu(\Omega_k)-\mu(\Omega_{k+1})\Big)
=(\rho-1) \sum_{k<0}\rho^{k}  \Big(\mu(\Omega_k)-\mu(\Omega_0)\Big) \]
Hence,
\[ \int_{\Omega_+}g^{2}d\mu 
\le  s \,\sup (f-m)_+^{2}+ \sum_{k<0} c (\rho-1)\rho^k (\mu(\Omega_k)-s).\]
In order to use our hypothesis, note that it implies that for every $A$ 
with measure at most $1/2$, one has
$\capa_\mu(A) \ge (\mu(A)-s)/\gamma(s).$ Indeed this is obvious if $s\ge \mu(A)$,
whereas if $s\le \mu(A)$, $\capa_\mu(A) \ge \mu(A)/\gamma(\mu(A)) \ge (\mu(A)-s)/
 \gamma(s)$ by the monotonicity of $\gamma$. 
Thus choosing
$$ g_k=
 \min\left(1,\left( \frac{g-\sqrt{c\rho^{k-1}}}
{\sqrt{c\rho^{k}}-\sqrt{c\rho^{k-1}}}  \right)_+  \right),$$
we have
\begin{eqnarray*}
\mu(\Omega_k)-s 
& \le & 
\gamma(s) \capa_\mu(\Omega_k) \le \gamma(s) \int |\nabla g_k|^{2}d\mu \\
& \le &
\gamma(s) \int_{\Omega_{k-1} \setminus \Omega_k}
\frac{|\nabla g|^{2}}{c\rho^{k-1} (\sqrt\rho-1)^{2}}d\mu.
\end{eqnarray*}
Summing upon $k<0$ we obtain
\begin{eqnarray*}
   \int_{\Omega_+}g^{2}d\mu 
  &\le & s \,\sup (f-m)_+^{2}+\gamma(s)\frac{\rho(\rho-1)}{(\sqrt\rho-1)^{2}}
  \sum_{k<0}  \int_{\Omega_{k-1}\setminus \Omega_k} |\nabla g|^{2}d\mu\\
  &\le &\gamma(s)\rho\frac{\sqrt\rho+1}{\sqrt\rho-1}\int_{\Omega_+} |\nabla f|^{2}d\mu
 + s \,\sup (f-m)_+^{2}.
\end{eqnarray*}
Summing up with a similar estimate for $\Omega_-$ and optimizing on $\rho$ gives
a slightly better estimate than the claimed one.
\end{proof}

\begin{theorem}
   Let $\mu$ be a probability measure on $\dR$. Assume that it is absolutely
continuous with respect to Lebesgue measure 
and denote by $\rho_\mu$ its density. Let $m$ be a median of $\mu$.
Let $\beta:(0,1/2)\to \dR^{+}$ be non-increasing.
Let $C$ be the optimal constant such that for all $f:\dR\to \dR$ and $s\in(0,1/4)$,
$$\var_\mu(f) \le C\beta(s) \int |\nabla f|^{2}d\mu +s\, \osc(f)^{2}.$$
Then $\frac14 \max(b_-,b_+)  \le  C \le 12 \max (B_-,B_+)$, where
\begin{eqnarray*}
   b_+&=&\sup_{x>m} \mu([x,+\infty)) \frac{1}{\beta(\mu([x,+\infty))/4)} 
              \int_m^{x}\frac{1}{\rho_\mu} \\
  b_-&=&\sup_{x<m} \mu((-\infty,x]) \frac{1}{\beta(\mu((-\infty,x])/4)} 
              \int_x^{m}\frac{1}{\rho_\mu} \\
    B_+&=&\sup_{x>m} \mu([x,+\infty)) \frac{1}{\beta(\mu([x,+\infty)))} 
              \int_m^{x}\frac{1}{\rho_\mu} \\
        B_-&=&\sup_{x<m} \mu((-\infty,x]) \frac{1}{\beta(\mu((-\infty,x]))} 
              \int_x^{m}\frac{1}{\rho_\mu} \cdot
\end{eqnarray*}
\end{theorem}
\begin{proof}
We start with the lower bound on $C$. We have seen that the weak spectral gap 
inequality ensures that for all $\Omega$ with $\mu(A)\le 1/2$ and 
$A\subset \Omega$,
one has $\capa_\mu(A,\Omega) \ge \mu(A)/(4C\beta(\mu(A)/4))$. 
Let $x>m$ and apply this
inequality with $A=[x,+\infty)$ and $\Omega=(m,+\infty)$.
It is easy to check 
that 
$\capa_\mu([x,+\infty),(m,+\infty))=1/\int_m^{x}1/\rho_\mu$.
This yields $C\ge b_+/4$. A similar argument on the other side of the median $m$
also gives  $C\ge b_-/4$.

   For the upper bound, we follow the argument of the proof of Theorem~\ref{crit-capa}
with some modification. We start with writing that
$$ \var_\mu(f) \le \int_m^{+\infty} |f-f(m)|^{2}d\mu 
    +\int_{-\infty}^{m}|f-f(m)|^{2}d\mu.$$
We work separately on the right and on the left of $m$. We explain only for
the right side; the left one is similar.
To proceed the argument in the same way  we need to check that any 
$A\subset (m,+\infty)$ verifies 
$$ \capa_\mu(A,(m,+\infty))\ge \frac{\mu(A)}{B_+\,\beta(\mu(A))}.$$
By hypothesis the above inequality holds  when $A=[x,+\infty)$. It follows
that it is valid for general $A$. Indeed,
for  any $A\subset (m,+\infty)$ one has $\capa_\mu(A,(m,+\infty)) =
\capa_\mu([\inf A,+\infty),(m,+\infty))$.
Since $\mu(A)\le \mu([\inf A,+\infty))$ and  $t\mapsto t/\beta(t)$
is non-decreasing the above inequality for half-lines implies it for general sets.
\end{proof}

\begin{corollary}\label{cor:wp}
   Let $d\mu(x)=e^{-\Phi(x)}dx,$ $x\in \dR$ be a probability measure.
 Let $\varepsilon \in (0,1)$. Assume that
there exists an interval $I=(x_0,x_1)$ containing a median $m$ of $\mu$ such 
that $|\Phi|$ is bounded on $I$, and $\Phi$ is twice differentiable outside $I$
with 
$$\Phi'(x)\neq 0 \quad \mathrm{and} \quad  
\frac{|\Phi''(x)|}{\Phi'(x)^{2}}\le 1-\varepsilon, \qquad x\not \in I.$$
Let  $\beta$ be a decreasing function on $(0,1/2)$.
 Assume that there exists $c>0$ such that 
for all $x\not\in I$ one has 
$$\beta\left( \frac{e^{-\Phi(x)}}{\varepsilon |\Phi'(x)|}\right)
  \ge \frac{c}{\Phi'(x)^{2}}.$$
Then $\mu$ satisfies a weak Poincar\'e inequality with function 
$C\beta$ for some constant $C>0$.
\end{corollary}

\begin{proof}
   We evaluate the quantity $B_+$ in the above theorem. The study of $B_-$ is
similar. For $x\ge x_1$, we have 
$$ \left( \frac{e^{\Phi}}{\Phi'}\right)'(x)
   =e^{\Phi(x)}\left(1-\frac{\Phi''(x)}{\Phi'(x)^{2}}\right)
   \ge \varepsilon e^{\Phi(x)} .$$
Therefore by integration
$$ \int_m^{x}e^{\Phi} \le  \int_m^{x_1}e^{\Phi} + \int_{x_1}^{x}e^{\Phi}
             \le (x_1-m) e^{M}  + \frac{1}{\varepsilon}
         \left(\frac{e^{\Phi(x)}}{\Phi'(x)}-\frac{e^{\Phi(x_1)}}{\Phi'(x_1)} \right),$$
where $M=\sup\{|\Phi(x)|;\; x\in I\}$. Similar calculations give
$$            (2-\varepsilon)      e^{-\Phi(x)}         \ge
\left(- \frac{e^{-\Phi}}{\Phi'}\right)'(x)
   \ge \varepsilon e^{-\Phi(x)} .$$
Note  that $\lim_{+\infty} e^{-\Phi}/\Phi'=0$. Indeed 
this quantity is positive, since $\Phi'$ cannot change sign, and decreasing
by the above bound. The limit has to be zero otherwise $e^{-\Phi(x)}$ would 
behave as $c/x$ and would not be integrable. We obtain by integration for $x\ge x_1$,
$$ \mu([x,+\infty))\le \frac{e^{-\Phi(x)}}{\varepsilon \Phi'(x)}
\le \frac{2-\varepsilon}{\varepsilon}\mu([x,+\infty)).$$
Combining these bounds on $\int_m^{x}e^{\Phi}$ and $\mu([x,+\infty))$ it is not hard
to show that $B_+$ is finite.
 \end{proof}
 
\begin{example}
For $\alpha>0$, the measure $dm_\alpha(t)= \alpha(1+|t|)^{-1-\alpha}dt/2,\, t\in \dR$
satisfies the weak spectral gap inequality with $\beta(s)=c_\alpha s^{-2/\alpha}$.
This was proved differently in \cite{rockw01wpil}, our next result improve on theirs.
\end{example}
\begin{example}
For $p\in(0,1)$, the measure $d\nu_p(t)=e^{-|t|^{p}}/(2\Gamma(1+1/p)),\; t\in \dR$
satisfies the inequality with $\beta(s)=d_p \log(2/s)^{\frac{2}{p}-2}.$
\end{example}
\begin{remark} 
In the above examples, the functions $\beta$ are best possible up to
 a multiplicative constant (we could write an analogue of the previous 
corollary, providing a necessary condition for a weak Poincar\'e inequality
to hold with $\beta$, with a similar proof). Since these functions $\beta$ 
satisfy the doubling condition,  our  theorem  describes all real measures
enjoying the same functional inequality.
\end{remark}

\section{Tensorisation} It is classical that the Poincar\'e inequality enjoys
the tensorisation property. When $\beta$ has infinite limit at 0, the weak spectral 
gap inequality does not tensorise. We shall give geometric evidence for this
in the section related to isoperimetry. However if $\mu$ satisfies the inequality
with a function $\beta$, then $\mu^{n}$  satisfies a weak spectral gap inequality
  with a worse function.

\begin{theorem}
   Assume that for every $f:M\to \dR$ and every $s\in (0,1/4)$ one has
$$ \var_\mu(f) \le \beta(s) \int |\nabla f|^{2}d\mu + s \,\osc(f)^{2}.$$
Let $n\ge1$. Then for every $f:M^{n}\to \dR$ and every $s\in (0,1/4)$ one has
$$ \var_{\mu^{n}}(f) 
  \le \beta\Big(\frac{s}{n}\Big) \int |\nabla f|^{2}d\mu^{n} + s \,\osc(f)^{2}.$$
\end{theorem}
\begin{proof}
By the sub-additivity property of the variance,
$$ \var_{\mu^{n}}(f) \le \sum_{i=1 }^{n} \int \var_\mu\Big(y_i\mapsto
 f(x_1,\ldots, x_{i-1},y_i,x_{i+1},\ldots,x_n)\Big) \prod_{j\neq i}d\mu(x_j).$$  
For each $i$ the inner variance is at most 
$$ \beta(s) \int |\nabla_i f|^{2}(x_1,\ldots,y_i,\ldots,x_n) d\mu(y_i)+s \,
 \osc\Big(y_i\mapsto f(x_1,\ldots,y_i,\ldots,x_n)\Big)^{2}.$$
The latter oscillation is less than or equal to $\osc(f)$. Summing up we arrive at
$$\var_{\mu^{n}}(f) \le \beta(s) \int |\nabla f|^{2}d\mu^{n}+ns \, \osc(f)^{2},$$
for all $s\in (0,1/4).$
\end{proof}

\section{Isoperimetric inequalities}
For $h>0$ we denote the $h$-enlargement of a set $A\subset M$ in the 
geodesic distance by $A_h$.
The boundary measure in the sense of $\mu$ is by definition
$$ \mu_s(\partial A)=\liminf_{h\to 0}\frac{\mu(A_h\setminus A)}{h}.$$
The isoperimetric function encodes the minimal boundary measure of sets
of prescribed measures:
$$I_\mu(a)=\inf\{\mu_s(\partial A);\; \mu(A)=a\},\qquad a\in [0,1].$$

It was shown by R\"ockner and Wang that in the diffusion case,
a weak spectral gap inequality for $\mu$ implies an isoperimetric 
inequality. We state here a consequence of their results.

\begin{theorem}[\cite{rockw01wpil}]
   Let $\mu$ be a probability measure on $(M,g)$, with density $e^{-V}$ with respect
to the volume measure. Assume that  $V$ is $C^{2}$ and such 
that $Ricci+ \nabla\nabla V\ge Rg$
for some  $R\le 0$. If  $\mu$ satisfies a weak spectral gap inequality
with function $\beta$, with $\beta(1/8)\ge\varepsilon>0$,
 then for every measurable $A\subset M$,
$$ \mu_s(\partial A) \ge c(\varepsilon,R) \frac{p}{\beta(p/2)},$$
where $p=\mu(A)(1-\mu(A))\ge \min(\mu(A),\mu(A^{c}))/2$.
\end{theorem}

\begin{remark} 
Comparing with a result of R\"ockner and Wang, showing that an isoperimetric
inequality implies a weak spectral gap inequality,
 one notices that $\sqrt\beta$ is expected
in the denominator (in the method, this loss comes from the necessity 
to estimate the underlying semi-group for large time instead of small time).
\end{remark}

\begin{corollary}
   Under the hypothesis of the above theorem, the following isoperimetric inequality
 holds for all $n\ge 1$. For all $A\subset M^{n}$, one has
$$ \mu^{n}_s(\partial A) \ge c(\varepsilon,R) \frac{p}{\beta(p/(2n))},$$
where $p=\mu^{n}(A)(1-\mu^{n}(A))$. 
\end{corollary}
\begin{proof}
   The tensorisation result of the previous section provides a weak spectral
  gap inequality for $\mu^{n}$ with function $\beta(s/n)$. The latter
  theorem then applies. Note that  the differential
  hypothesis on the density of $\mu$ remains valid for $\mu^{n}$. We also used
  $\beta(1/(8n))\ge \beta(1/8)\ge \varepsilon$. 
\end{proof}
In the non-trivial cases when $\lim_0\beta=+\infty$ the above lower bound
of $I_{\mu^{n}}$ tends to zero as $n$ increases. This has to be, as  
the following  consideration of product sets shows. We shall assume that
$I_\mu(t)=I_\mu(1-t)$ for all $t$ (this is very natural, since regular 
sets have the same boundary measure as their complement).
First note that for all $n\ge 1$, $h>0$ and $A\subset M$
one has $(A^{n})_h\subset (A_h)^{n}$, where $A^{n}\subset M^{n}$ is the
cartesian product of $n$ copies of $A$. Combining this with the definition
of the boundary measure yields
$$\mu^{n}_s(\partial (A^{n})) \le n\mu(A)^{n-1}\mu_s(\partial A).$$
Taking infimum on $A$ with prescribed measure, we get $I_{\mu^{n}}(a^{n})\le n a^{n-1}
I_\mu(a)$ for all $a\in (0,1)$.
Thus for any fixed $t\in(0,1)$ one has when $n\ge \log(1/t)/\log(2)$
\begin{eqnarray*}
   I_{\mu^{n}}(t)&\le& n t^{1-\frac1n}I_\mu(t^{\frac1n}) \le 2nt I_\mu(1-t^{\frac1n})\\
    &=& 2nt\, I_\mu\left(\frac{\log(1/t)}{n}(1+\varepsilon_t(n)) \right)\\
    &=& 2t \log\Big(\frac1t\Big)\, \Theta\left(\frac{\log(1/t)}{n}(1+\varepsilon_t(n))
 \right)
 (1+\varepsilon_t(n)),
\end{eqnarray*}
where $\lim_n \varepsilon_t(n)=0$ and  $I_\mu(u)=u\Theta(u)$.
If $\Theta$ tends to zero at zero then $\lim_n I_{\mu^{n}}(t)=0$ with corresponding
speed. 

\medskip
 For even measures on $\dR$ with positive density on
a segment, Bobkov and Houdr\'e \cite[Corollary 13.10]{bobkh97cbis}
 proved that solutions to the isoperimetric 
problem can be found among half-lines, symmetric segments and 
their complements. More precisely, if $\rho_\mu$ is the density and $R_\mu$ 
the distribution function of $\mu$, then denoting $J_\mu=\rho_\mu \circ R_\mu^{-1}$,
one has for $t\in(0,1)$
$$ I_\mu(t)=\min\left(J_\mu(t),2J_\mu\Big(\frac{\min(t,1-t)}{2} \Big) \right).$$
This readily applies to our previous examples.

\begin{example}
For the measures $dm_\alpha(t)=\alpha(1+|t|)^{-1-\alpha}/2$ one gets
$J_{m_\alpha}(t)=\alpha 2^{1/\alpha}\min(t,1-t)^{1+1/\alpha}$, and thus for
 $t\in(0,1/2),$ 
$$I_{m_\alpha}(t)=\alpha t^{1+1/\alpha}.$$
The results of this section do not apply to $m_\alpha$ for lack of regularity.
However for an even unimodal smoothed perturbation $\tilde{m}_\alpha$, up to 
a numerical constant, the
same isoperimetric and weak spectral gap inequality hold. So there are constants
such that for  $t\le 1/2$ and $n\ge\log(1/t)/\log2$ one has
$$c_1(\alpha)\, t \left(\frac{t}{n} \right)^{2/\alpha}\le I_{\tilde{m}^{n}_\alpha}(t)
 \le c_2(\alpha)\, t \frac{\log(1/t)^{1+1/\alpha}}{n^{1/\alpha}}.$$
\end{example}

\begin{example}
For $p\in(0,1)$, and $d\mu_p(t)=\exp(-|t|^{p})/(2\Gamma(1+1/p))$ similar estimates
can be done. For $t\le 1/2$, $I_{\nu_p}(t)$ is comparable to $t (\log(1/t))^{1-1/p}$.
So for a suitable smoothed version of this measure, one gets 
$$d_1(p)\, t\left(\log\Big(\frac{n}{t}\Big) \right)^{2(1-1/p)}\le I_{\tilde\nu_p^{n}}(t)
\le d_2(p)\, t\log(1/t) \left(\log\Big(\frac{n}{\log(1/t)}\Big) \right)^{1-1/p},$$
which guarantees a convergence to zero with logarithmic speed in the dimension. 
\end{example}

\section{Concentration of measure}
In this section, we shall derive concentration inequalities, that is lower bounds
on the measure of enlargements of rather large sets, or equivalently
deviation inequalities for Lipschitz functions. They can be approached via
isoperimetric inequalities, which quantify the measure of infinitesimal enlargements.
In our setting, we have seen in the previous section that the available methods
provide loose isoperimetric bounds. Hence we come back to simpler and 
more robust techniques.
It is known, since Gromov and Milman \cite{gromm83taii}, that  a Poincar\'e inequality
yields exponential concentration. See e.g. \cite{ledoCMLS} for subsequent 
developments. We show how a weak spectral gap inequality 
can be used to derive deviation inequalities for Lipschitz functions. 
Among the various available methods used for Poincar\'e inequalities, the one
in Aida, Masuda and Shigekawa  \cite{aidams94lsie} is the most adapted.

\begin{theorem}
   Let $\mu$ satisfy a weak spectral gap inequality with function $\beta$.
 Let $F:M\to \dR$ be a $L$-Lipschitz function with median $m$. 
 Then for $k\ge 1$ and $s\in(0,1/4)$, one has
 \begin{equation}\label{eq:rec}
    \mu(F-m>k) 
 \le \frac{s}{1+L^{2}\beta(s)}
                  + \mu(F-m>k-1)  \left(1-\frac{1}{2(1+L^{2}\beta(s))}\right).
 \end{equation}
Consequently
\begin{equation}\label{eq:dev}
 \mu(F-m>k) \le 2s+\frac{\sqrt e}{2} \exp\Big(\frac{-k}{4L\sqrt{\beta(s)}}\Big).
\end{equation}
Thus
$ \mu(|F-m|>k)\le 6\,\Theta(k/L),$
where $$ \Theta(u)=\inf\Big\{s\in(0,1/4];\;
  \exp\Big(\frac{-u}{4\sqrt{\beta(s)}}\Big)\le s\Big\}$$
tends to $0$ when $u$ tends to infinity.
\end{theorem}
\begin{proof}
For notational convenience assume that $m=0$.
   Let $\varepsilon>0$. Let $\Phi:\dR\to \dR^{+}$ be a non-decreasing smooth 
function with $\Phi_{|(-\infty,\varepsilon]}=0$, $\Phi_{|[1-\varepsilon,+\infty)}=1$
and $\|\Phi'\|_\infty\le 1+3\varepsilon$. Set $\Phi_k(t)=\Phi(t-k+1)$. We apply the
weak Poincar\'e inequality  to $\Phi_k(F)$. Since $\mathbf 1_{(k-1,+\infty)} \ge \Phi_k 
\ge \mathbf 1_{[k,+\infty)}$ one has 
$$ \int \Phi_k(F)^{2}d\mu\ge \mu(F\ge k),\qquad
 \left(\int \Phi_k(F) d\mu\right)^{2} \le \mu(F>k-1)^{2}.$$
Almost surely one has $|\nabla \Phi_k(F)| \le |\Phi_k'(F)|\cdot |\nabla F| \le
(1+3\varepsilon)L\mathbf1_{k-1<F<k}$.
Therefore, letting $\varepsilon$ to zero, the inequality
$$\var(\Phi_k(F))\le \beta(s)\int |\nabla \Phi_k(F)|^{2}d\mu +s \osc(\Phi_k(F))^{2},$$
readily implies
$$ \mu(F>k)-\mu(F>k-1)^{2}\le L^{2}\beta(s) \Big(\mu(F>k-1)-\mu(F>k) \Big)+s.$$
Rearranging 
$$  \mu(F>k) \le \frac{s}{1+L^{2}\beta(s)}
                  + \mu(F>k-1) \frac{\mu(F>k-1)+L^{2}\beta(s)}{1+L^{2}\beta(s)}.$$
The first claimed inequality follows from the above and 
$\mu(F>k-1)\le \mu(F>0) \le 1/2$. 
Iterating this inequality $k$ times gives
\begin{eqnarray*}
   \mu(F>k)& \le& \frac{s}{1+L^{2}\beta(s)} 
\left[ 1+\left(1- \frac{1}{2(1+L^{2}\beta(s))}\right) 
           + \cdots + \right. \\
   &&  \left. \Big(1-\frac{1}{2(1+L^{2}\beta(s))} \Big)^{k-1} \right]
   + \left(1-\frac{1}{2(1+L^{2}\beta(s))} \right)^{k}\mu(F>0)\\
  &\le&\frac{s}{1+L^{2}\beta(s)}\cdot\frac{1}{1-(1-\frac{1}{2(1+L^{2}\beta(s))})}
   + \frac12  \left(1-\frac{1}{2(1+L^{2}\beta(s))} \right)^{k}\\
   &\le& 2s+\frac12 \exp\left(\frac{-k}{2(1+L^{2}\beta(s))}\right).
\end{eqnarray*}
Note that this is also true when $k=0$.
Let $\lambda >0$ and apply the latter bound to the $\lambda L$-Lipschitz function
$\lambda F$ with median $0$. Denoting by $[x]$ the integer part of $x$, we get  
\begin{eqnarray*}
   \mu(F>k)&=& \mu(\lambda F>\lambda k) \le \mu(\lambda F > [\lambda k]) \\
     &\le& 2s 
      +\frac12 \exp\Big(\frac{-[\lambda k]}{2(1+\lambda^{2}L^{2} \beta(s))} \Big) \\
      &\le& 2s 
      +\frac12 \exp\Big(\frac{-\lambda k}{2(1+\lambda^{2}L^{2} \beta(s))}
            +\frac{1}{2(1+\lambda^{2}L^{2} \beta(s))} \Big) \\  
      &\le& 2s 
      +\frac{\sqrt e}{2}
       \exp\Big(\frac{-\lambda k}{2(1+\lambda^{2}L^{2} \beta(s))} \Big).
\end{eqnarray*}
Choosing $\lambda=1/(L\sqrt{\beta(s)})$ establishes (\ref{eq:dev}).
The rest of the statement  easily follows.
\end{proof}
Next we give a few examples. 
\begin{example}
If $\beta$ has a finite limit at $0$ then taking
$s=0$ in (\ref{eq:dev}) recovers the well known exponential deviation inequality. 
\end{example}

\begin{example} 
Let $F:\dR^{n}\to \dR$ be a 1-Lipschitz function with median $m$.
We consider on $\dR^{n}$ the n-fold product of $dm_\alpha(t)=\alpha(1+|t|)^{-1-\alpha}/2
\, dt$ denoted $m_\alpha^{n}$. Since this measure satisfies a weak spectral
gap inequality with $\beta(s)=c_\alpha (s/n)^{-2/\alpha}$,  the deviations
of $\mu(F-m>k)$ are controlled by 
$$\inf_{s\in(0,1/4)} 2s+\frac{\sqrt e}{2} \exp\Big(\frac{-k s^{1/\alpha}}
{4\sqrt{c_\alpha} n^{1/\alpha}} \Big).$$
Setting $t=k/(4\sqrt{c_\alpha} n^{1/\alpha})$, we choose 
$s=(\alpha \log(t)/t)^{\alpha}$. It is in the interval $(0,1/4)$ provided
$t$ is larger than a constant $t_1(\alpha)$. Under this hypothesis the infimum
is bounded from above by 
$$ 2(\alpha \log(t)/t)^{\alpha}+1/t^{\alpha}.$$
Therefore there exists constants $t_0(\alpha)>e$ and $C(\alpha)$ 
such that for $t\ge t_0(\alpha)$
\begin{equation}\label{eq:ma}
 m_\alpha^{n}(|F-m|>t n^{1/\alpha}) \le 
C(\alpha)  \left( \frac{ \log(t)}{t}\right)^{\alpha}.
\end{equation}
This is valid provided $4t\sqrt{c_\alpha}n^{1/\alpha}\in \mathbb N$ but extends 
to general values of $t$, with  slightly worse constants.
As we show next, this estimate is  correct up to the $\log$ factor. 
Presumably, this point  could be improved by optimizing
in $s$ the recursion formula (\ref{eq:rec}). 

Let us prove that (\ref{eq:ma}) is very close to the truth, by adapting
Talagrand's argument. It consists in analyzing product sets. First note that
if $A\subset \dR^{n}$ has measure at least $a\ge 1/2$ then $0$ is a median of the
distance function 
$x\mapsto d(x,A)$. Since the latter is 1-Lipschitz, (\ref{eq:ma}) applies
and gives,
\begin{equation}\label{eq:ma2}
m_\alpha^{n}(A_{tn^{1/\alpha}}) 
   \ge 1-C(\alpha) \left( \frac{\log t}{t}\right)^{\alpha}.
\end{equation}
We show that this is close to optimal by choosing a specific product set.
Namely we take $A=(-\infty,R^{-1}(a^{1/n})]^{n}$, where $R=R_{m_\alpha}$ is the 
distribution function of $m_\alpha$ and $R^{-1}$ is its reciprocal function.
By definition $m_\alpha^{n}(A)=a$. For $h>0$, its $h$-enlargement satisfies
$$
m_{\alpha}^{n}(A_h)
\le 
m_{\alpha}^{n}(A+[-h,h]^{n})
\!=\!
m_{\alpha}^{n}((-\infty,R^{-1}(a^{1/n})+h]^{n})
\!=\!
R\Big(R^{-1}(a^{1/n})+h\Big)^{n} \!\! .
$$
The function $R$ is explicitly computed. The latter estimate thus becomes
\begin{eqnarray*}
  m_{\alpha}^{n}(A_h)&\le&
     \left(1-\frac{1}{2\left(h+(2(1-a^{1/n}))^{-1/\alpha}
                       \right)^{\alpha}} 
    \right)^{n} \\
           &\le& \exp\left(\frac{-n}{2\left(h+(\frac{2}{n}\log(\frac{1}{a})
         +O(\frac{1}{n^{2}}))^{-1/\alpha}
                       \right)^{\alpha}}\right) \\
    &=& \exp\left(\frac{-1}{2\left(\frac{h}{n^{1/\alpha}}+(2\log(\frac{1}{a})
         +O(\frac{1}{n}))^{-1/\alpha}
                       \right)^{\alpha}}\right).
\end{eqnarray*}
We think of $A$ and $h$ as depending on $n$. The above bound shows that when
$n$ is large and $h<<n^{1/\alpha}$ the measure of  $A_h$ is essentially
equal to $a= m_{\alpha}^{n}(A)$. This confirms  that $h=tn^{1/\alpha}$ is the
right scale of enlargement. In this scale we  have
$$
  m_{\alpha}^{n}(A_{tn^{1/\alpha}})\le
      \exp\left(\frac{-1}{2\left(t+(2\log(\frac{1}{a})
         +O(\frac{1}{n}))^{-1/\alpha}
                       \right)^{\alpha}}\right)
  \le 1-\frac{c_\alpha}{t^{\alpha}},
$$
when $t\ge t_2(\alpha)$. Comparing this with Inequality (\ref{eq:ma2}) 
proves the tightness of our bounds.
\end{example}

\begin{example}
Finally, we consider the measures 
$\nu_p^{n}=\left(d_p e^{-|t|^{p}}dt\right)^{\otimes n}$, for $p\in(0,1)$. 
We have shown that they satisfy the weak Poincar\'e inequality with $\beta(s)=k_p 
\log(2n /s)^{(2/p)-2}.$ Therefore the deviations of Lipschitz functions are
controlled by 
$$\inf_{s\in(0,1/4)} 2s+\exp\left(\frac{-k (\log(2n/s))^{1-1/p}}{4\sqrt{k_p}}  \right).$$
We look for a value of $s$ such that the two terms are of similar size. We are
inspired by the case $p=1/2$ where explicit calculations can be done.

If $k\ge (\log n)^{1/p}$ we set $s=2e^{-k^{p}}$. The above infimum is at most
(denoting by $c_p$ a quantity depending only on $p$ and that may be different
 in different occurrences)
 \begin{eqnarray*}
  \nu_p^{n}(F-m>k) &\le&   4e^{-k^{p}}+e^{-kc_p (\log n+k^{p})^{1-1/p}}\\
 &\le&  4e^{-k^{p}}+e^{-kc_p (2k^{p})^{1-1/p}}\\
 &\le& 5 e^{-c_p k^{p}}.
 \end{eqnarray*}
 Here
we did not check that the chosen $s$ is less than $1/4$, since otherwise the bound
is trivial.

If $k\le (\log n)^{1/p}$ we set $s=2e^{-k (\log n)^{1-1/p}}$. We get
\begin{eqnarray*}
  \nu_p^{n}(F-m>k) &\le& 4 e^{-k (\log n)^{1-1/p}}
           + e^{-kc_p (\log n+k(\log n)^{1-1/p})^{1-1/p}}\\
   &\le&  4 e^{-k (\log n)^{1-1/p}}
           + e^{-kc_p (2 \log n)^{1-1/p}}\\
    &\le&  5 e^{-c_p k (\log n)^{1-1/p}}.
\end{eqnarray*}
As a conclusion we obtained 
$$ \nu_p^{n}(|F-m|>k)\le 10 
\exp\left(\frac{-c_p k}{\max(k^{p},\log n)^{\frac1p-1}}\right).$$
In particular, for $\varepsilon$ fixed and $n$ large, it is enough to take
$k\ge c_p (\log\frac{10}{\varepsilon})(\log n)^{\frac1p-1}$ in order to ensure
$   \nu_p^{n}(|F-m|>k)  \le \varepsilon$. 
\end{example}

\begin{remark}
   Theorem~2.4 of \cite{wangz03wpid}  also derives concentration inequalities
from a weak spectral gap inequality, but they are different from ours.
 Comparing their Corollary~2.5 with the above examples shows that our
 result is  sharper. The main technical 
reason for this is that the final step of our proof (which reintroduces homogeneity,
as it was destroyed by the cut-off method) is not performed.
 Combining their method and the optimization on a scaling factor $\lambda$ provides
a slightly better estimate than ours. Let $c\in(0,1/2)$, then with the notation
of the theorem
$$\left(\frac12-c\right)^{2} \frac{k^{2}}{4}\le \left[
 \log\left(\frac{1}{2\mu(F-m>k)}\right)+\frac12-c \right]
       \int_{\mu(F-m>k)}^{\frac12} \frac{\beta(cs)}{s}ds. $$
 In general though, the integral can only be
estimated by $$\beta\big(c\mu(F-m>k)\big)\log\left(\frac{1}{2\mu(F-m>k)}\right).$$
 This recovers our bound. For the measures $m_\alpha$
the integral can be computed and one gets a better decay, by a different
power on the $\log$-term. In the case of $\nu_p$ the explicit computation
does not improve on our result. 
\end{remark}

\section{Concave potentials of power type}

In this section we apply our methods to products of probability measures 
on $\dR$, $d\mu_\Phi(x) = Z_\Phi^{-1} e^{-\Phi(|x|)}dx$, where
$\Phi$ satisfies the following assumption:

\begin{assumption}[H]
\begin{description}
\item{$(i)$} $\Phi : \R^+ \rightarrow \R^+$ is an increasing concave
   function with $\Phi(0)=0$ and ${\cal C}^2$ in a neighborhood of $+
   \infty$.
\item{$(ii)$} There exists $B>1$ such that for $x$ large enough
   $\Phi(2x) \geq B \Phi(x)$.
\item{$(iii)$} There exists $C>0$ such that for $x$ large enough $|x
   \Phi''(x)| \leq C \Phi'(x)$.
\end{description}
\end{assumption}

Hypothesis $(H)$ naturally generalizes the power potentials
$\Phi_p(x)=|x|^p$, $p \in (0,1)$. In particular it is not hard to
check that $\Phi_{p,\beta} = |x|^p \log(\gamma + |x|)^\alpha$ with $p
\in (0,1)$, $\alpha >0$ and $\gamma=e^{2\alpha/(1-p)}$ verifies
Hypothesis $(H)$ with $B=2^p$ and $C=1$.

\begin{remark} \label{rem:H}
   Assertion $(ii)$ of $(H)$ yields $\lim_{+\infty} \Phi = +\infty$ and by induction
   for large $x$
   \begin{equation}
      \label{eq:2}
      2\Phi(x)\le \Phi(B'x),
   \end{equation}
 with $B'=2^{1+\log2/\log B}>1$.
   On the other hand, since $\Phi$ is concave and $\Phi(0)=0$, $(ii)$
   also implies that 
\begin{equation}\label{eq:1}
(B-1) \Phi(x) \leq \Phi(2x) - \Phi(x) \leq \int_x^{2x} \Phi'
\leq x \Phi'(x) \leq \int_0^x \Phi' = \Phi(x) 
\end{equation}
where the left inequality is valid for $x$ large enough, and the other ones for 
$x\ge0$ (when $\Phi$ is not differentiable, $\Phi'(x)$ stands for the right derivative).
Together with $(iii)$ this result implies that $\displaystyle
\frac{|\Phi''(x)|}{\Phi'(x)^2} \leq \frac{C}{x \Phi'(x)} \leq
\frac{C}{(B-1)\Phi(x)}.  $ Hence, $\lim_{+ \infty}
\frac{|\Phi''|}{(\Phi')^2}= 0$. Also, combining the concavity assumption with 
\eqref{eq:2} and \eqref{eq:1} yields for $x$ large enough
\begin{equation}
   \label{eq:3}
   \Phi'(x)\ge \Phi'(B'x)\ge B'' \Phi'(x),
\end{equation}
where $B''\in(0,1)$ depends only on $B$.
\end{remark}

Now we prove that $\mu_\Phi$ satisfies a weak Poincar\'e inequality
with appropriate function $\beta$.

\begin{proposition}\label{prop:wp}
   Let $d\mu_\Phi(x) = Z_\Phi^{-1} e^{-\Phi(|x|)}dx$ be a probability
   measure on $\dR$.  Assume that $\Phi$ verifies Hypothesis $(H)$.
   Then there exists a constant $c_\Phi>0$ such that $\mu_\Phi$
   satisfies a weak Poincar\'e inequality with function $c_\Phi \beta$
   where
   $$
   \beta(s) = \frac{1}{[\Phi' \circ \Phi^{-1}(\log \frac
     1s)]^2} ,\quad s\in(0,1/4).
   $$
\end{proposition}

\begin{proof}
   We use Corollary \ref{cor:wp}. From Hypothesis $(H)$ and the above
   remark there exists $A>0$ such that, for $x>A$
   $$
   \Phi'(x) \neq 0 \quad \mathrm{and} \quad
   \frac{|\Phi''(x)|}{\Phi'(x)^{2}}\le \frac 12 .
   $$
   Thus, we only have to check that $\beta
   \left(\frac{2e^{-\Phi(x)}}{|\Phi'(x)|} \right) \geq
   \frac{c}{\Phi'(x)^2}$ for some constant $c >0$ and $|x|$ large
   enough.
   
   It follows from Remark \ref{rem:H} that for $x$ large
   $\frac{\log(\Phi'(x))}{\Phi(x)} \leq
   \frac{\log(x\Phi'(x))}{\Phi(x)} \leq
   \frac{\log(\Phi(x))}{\Phi(x)}$.  Since $\lim_{+\infty} \Phi =
   +\infty$ we can deduce that $\lim_{+\infty} \frac{\log \Phi'}{\Phi}
   = 0$. Hence, for $x$ large enough one has
   $$
   \log \frac 12 + \log \Phi'(x) + \Phi(x) \geq \frac 12 \Phi(x).
   $$
   Now Equation \eqref{eq:2} implies that
   $$
   \Phi^{-1}(\log \frac 12 + \log \Phi'(x) + \Phi(x)) \geq
   \Phi^{-1}( \frac 12 \Phi(x)) \geq \frac{ x}{B'}.
   $$
 Since $\Phi'$ is non-increasing the above inequality and \eqref{eq:3} lead to
   \begin{eqnarray*}
   \beta \left(\frac{2e^{-\Phi(x)}}{ |\Phi'(x)|} \right) &=& \frac{1}
   {(\Phi')^2 \circ \Phi^{-1}(\log \frac 12 + \log \Phi'(x) +
     \Phi(x))}\\
    &\geq&\frac{1}{(\Phi')^{2}(x/B')} \geq  \left(\frac{B''}{\Phi'(x)}\right)^{2}
   \end{eqnarray*}
   for $x$ large enough. This achieves the proof.
\end{proof}

\begin{example}
   This result recovers the case $\Phi_{p}=|x|^p$, $p \in (0,1)$.  For
   $\Phi_{p,\alpha}=|x|^p \log(\gamma + |x|)^\alpha$ with $p \in
   (0,1)$, $\alpha >0$ and $\gamma=e^{2\alpha/(1-p)}$, one can easily
   see that $\mu_{p,\alpha}$ satisfies a weak Poincar\'e inequality
   with function asymptotically (when $s$ is small) behaving like
   $$
   \beta_{p,\alpha}(s) = \frac{1}{ (\log \frac
     1s)^{2(1-\frac 1p)} (\log \log \frac 1s)^\frac{2\alpha}{p}} .
   $$
\end{example}

We obtain the following concentration inequalities for $\mu_\Phi^{n}$:
\begin{proposition}
   Let $d\mu_\Phi(x) = Z_\Phi^{-1} e^{-\Phi(|x|)}dx$ be a probability
   measure on $\dR$ which verifies Hypothesis $(H)$. Then there exist
   $c_\Phi,\tilde{c}_\Phi,k_\Phi>0$ such that for any $n\ge 1$, any 1-Lipschitz
   function $F:\dR^{n}\to\dR$ and any integer $k\ge k_\Phi$ one has
\begin{eqnarray*}
\mu_\Phi^{n}(|F-m|>k)
& \le  & 6 
\exp \left(-c_\Phi k  \Phi' \circ \Phi^{-1} (\max(\Phi(k) , 2\log n))\right) \\
& \le & 6
\max\left(
e^{-\tilde{c}_\Phi \Phi(k)} , e^{-c_\Phi  k \Phi' \circ \Phi^{-1} (2\log n)}
\right) .
\end{eqnarray*}
where $m$ is a median of $F$.
\end{proposition}

\begin{proof}
   As in
   the previous section,
   since $\mu^n$ satisfies a weak Poincar\'e inequality with function
   $\beta(s)=c_\Phi/[\Phi' \circ \Phi^{-1}(\log \frac ns)]^2$, 
   the deviations of 1-Lipschitz functions are controlled by
   $$
   \inf_{s \in (0,1/4)} 2s + \frac{\sqrt e}{2} \exp\left(
      -\frac{4}{\sqrt{c_\Phi}} k \Phi' \circ \Phi^{-1}(\log \frac ns)
   \right) .
   $$
   Next we look for a value of $s$
   such that the two terms are of similar size.  We will denote by
   $c_\Phi$ a quantity depending only on $\Phi$ that may change from
   line to line. We work with $k$ large enough in order to be able
   to use the doubling condition in the following arguments.
   
   If $k \geq \Phi^{-1}(\log n)$ we set $s=e^{-\Phi(k)}$.  The above
   infimum is at most
 \begin{eqnarray*}
  \mu_\Phi^{n}(F-m>k) 
  &\le &  
  2e^{-\Phi(k)} + e^{- c_\Phi k\Phi' \circ \Phi^{-1}(\log n + \Phi(k))} \\
  &\le&  2e^{-\Phi(k)}+e^{- c_\Phi k \Phi'(k)}\\
  &\le&  3e^{- c_\Phi k \Phi'(k)} \le 3 e^{- \tilde{c}_\Phi \Phi(k)}.
 \end{eqnarray*}
 Here, we have used Equation \eqref{eq:2} 
 in order to get
 that 
 $$\Phi^{-1}(\log n + \Phi(k)) \leq \Phi^{-1}( 2 \Phi(k)) \leq B'k,$$
  and thus by \eqref{eq:3},  $\Phi' \circ
 \Phi^{-1}(\log n + \Phi(k)) \geq B''\Phi'(k)$.
  The last inequality comes from \eqref{eq:1}.

 If $k < \Phi^{-1}(\log n) $ we set $s=e^{-k \Phi' \circ
   \Phi^{-1}(\log n)}$.  Recall  first that   Inequality~\eqref{eq:1} asserts
that for $x\ge 0$ one has  $x \Phi'(x) \leq \Phi(x)$. 
Hence $ \Phi^{-1} (x) \Phi' \circ \Phi^{-1} (x) \leq x$
 and in turn it follows that 
$$k \Phi' \circ \Phi^{-1}(\log n) \leq
 \Phi^{-1}(\log n) \Phi' \circ \Phi^{-1}(\log n) \leq \log n.$$ 
 We get
\begin{eqnarray*}
  \mu_\Phi^{n}(F-m>k) 
  & \le & 
  2 e^{-k \Phi' \circ \Phi^{-1}(\log n)}
  + e^{-c_\Phi k \Phi' \circ\Phi^{-1}(\log n+k \Phi' \circ\Phi^{-1}(\log n))}\\
  & \le &  
  2 e^{-k \Phi' \circ \Phi^{-1}(\log n)}
  + e^{-c_\Phi k \Phi' \circ \Phi^{-1}(2 \log n)}\\
  & \le &  
  3 e^{-c_\Phi k \Phi' \circ \Phi^{-1}(2 \log n) }.
\end{eqnarray*}
The result easily follows.
\end{proof}

\medskip
\noindent
{\bf Acknowledgements:} We thank the referee for useful suggestions.

\begin{thebibliography}{10}

\bibitem{aida01egss}
S.~Aida.
\newblock An estimate of the gap of spectrum of {S}chr\"odinger operators which
  generate hyperbounded semigroups.
\newblock {\em J. Funct. Anal.}, 185(2):474--526, 2001.

\bibitem{aidams94lsie}
S.~Aida, T.~Masuda, and I.~Shigekawa.
\newblock Logarithmic {S}obolev inequalities and exponential integrability.
\newblock {\em J. Funct. Anal.}, 126(1):83--101, 1994.

\bibitem{bartcr04iibe}
F.~Barthe, P.~Cattiaux, and C.~Roberto.
\newblock Interpolated inequalities between exponential and {G}aussian,
  {O}rlicz hypercontractivity and application to isoperimetry.
\newblock {\em Revista Math. Iberoamericana}, To appear.

\bibitem{bartr03sipm}
F.~Barthe and C.~Roberto.
\newblock Sobolev inequalities for probability measures on the real line.
\newblock {\em Studia Math.}, 159(3):481--497, 2003.

\bibitem{bobkh00wdfc}
S.~G. Bobkov and C.~Houdr{\'e}.
\newblock Weak dimension-free concentration of measure.
\newblock {\em Bernoulli}, 6(4):621--632, 2000.

\bibitem{bobkh97cbis}
S.G. Bobkov and C.~Houdr{\'e}.
\newblock Some connections between isoperimetric and {S}obolev-type
  inequalities.
\newblock {\em Mem. Amer. Math. Soc.}, 129(616):viii+111, 1997.

\bibitem{borou83icnd}
A.~A. Borovkov and S.~A. Utev.
\newblock An inequality and a characterization of the normal distribution
  connected with it ({R}ussian).
\newblock {\em Teor. Veroyatnost. i Primenen}, 28(2):209--218, 1983.

\bibitem{gromm83taii}
M.~Gromov and V.~Milman.
\newblock A topological application of the isoperimetric inequality.
\newblock {\em Amer. J. Math.}, 105:843--854, 1983.

\bibitem{latao00bsp}
R.~Lata{\l}a and K.~Oleszkiewicz.
\newblock Between {S}obolev and {P}oincar\'e.
\newblock In {\em Geometric aspects of functional analysis}, number 1745 in
  Lecture Notes in Math., pages 147--168, Berlin, 2000. Springer.

\bibitem{ledoCMLS}
M.~Ledoux.
\newblock Concentration of measure and logarithmic {S}obolev inequalities.
\newblock In {\em S{\'e}minaire de Probabilit{\'e}s, XXXIII}, number 1709 in
  Lecture Notes in Math., pages 120--216, Berlin, 1999. Springer.

\bibitem{ledoCMP}
M.~Ledoux.
\newblock {\em The concentration of measure phenomenon}, volume~89 of {\em
  Mathematical Surveys and Monographs}.
\newblock American Mathematical Society, Providence, RI, 2001.

\bibitem{rockw01wpil}
M.~R{\"o}ckner and F.Y. Wang.
\newblock Weak {P}oincar{\'e} inequalities and ${L^{2}}$-conv\-erg\-ence 
rates of
  {M}arkov semigroups.
\newblock {\em J. Funct. Anal.}, 185:564--603, 2001.

\bibitem{schez00clb}
G.~Schechtman and J.~Zinn.
\newblock Concentration on the {$l\sp n\sb p$} ball.
\newblock In {\em Geometric aspects of functional analysis}, volume 1745 of
  {\em Lecture Notes in Math.}, pages 245--256. Springer, Berlin, 2000.

\bibitem{tala91niic}
M.~Talagrand.
\newblock A new isoperimetric inequality and the concentration of measure
  phenomenon.
\newblock In J.~Lindenstrauss and V.~D. Milman, editors, {\em Geometric Aspects
  of Functional Analysis}, number 1469 in Lecture Notes in Math., pages
  94--124, Berlin, 1991. Springer-Verlag.

\bibitem{wang03gbti}
F.-Y. Wang.
\newblock A generalized {B}eckner-type inequality.
\newblock {\em Preprint}.

\bibitem{wangz03wpid}
F.-Y. Wang and Q.~Zhang.
\newblock Weak {P}oincar{\'e} inequalities, decay of {M}arkov semigroups and
  concentration of measure.
\newblock {\em Preprint}.

\end{thebibliography}

\bigskip
\noindent
F. Barthe: Institut de Math\'ematiques. Laboratoire de Statistique et Probabilit\'es,
  UMR C 5583. Universit\'e Toulouse III. 118 route de Narbonne. 31062 Toulouse cedex 04.
  FRANCE.  

\noindent
Email: barthe@math.ups-tlse.fr

\medskip\noindent
P. Cattiaux: Ecole Polytechnique, CMAP, CNRS 756, 91128 Palaiseau Cedex FRANCE and
Universit\'e Paris X Nanterre, Equipe MODAL'X, UFR SEGMI, 200 avenue de la 
R\'epublique, 92001 Nanterre cedex, FRANCE.

\noindent
Email: cattiaux@cmapx.polytechnique.fr

\medskip\noindent
C. Roberto: Laboratoire d'analyse et math\'ematiques appliqu\'ees, UMR 8050. 
Universit\'es de Marne-la-Vall\'ee et de Paris 12  Val-de-Marne. 
Boulevard Des\-cartes, Cit\'e
Descartes, Champs sur Marne. 77454 Marne-la-Vall\'ee cedex 2. FRANCE

\noindent
Email: roberto@univ-mlv.fr

\end{document}